\documentclass[11pt]{article}
\usepackage{theorem}
\usepackage{amssymb}
\usepackage{amsmath}
\usepackage{amsfonts}

\def \RR {\mathbb R}
\def \EE {\mathbb E}

\def \E {\mathcal E}

\newtheorem{theorem}{Theorem}[section]
\newtheorem{lemma}[theorem]{Lemma}

\newtheorem{question}[theorem]{Question}
\newtheorem{proposition}[theorem]{Proposition}
\newtheorem{corollary}[theorem]{Corollary}
 {\theorembodyfont{\rmfamily}}

\begin{document}

\title{An isomorphic version of the slicing problem}
\author{B. Klartag\thanks{Supported by the Israel Clore
Foundation.}\\ School of Mathematical Sciences \\ Tel Aviv
University \\ Tel Aviv 69978, Israel}
\date{}
\maketitle

 \abstract{Here we show that any centrally symmetric
convex body $K \subset \RR^n$ has a perturbation $T \subset \RR^n$
which is convex and centrally symmetric, such that the isotropic
constant of $T$ is universally bounded. $T$ is close to $K$ in the
sense that the Banach-Mazur distance between $T$ and $K$ is
$O(\log n)$. If $K$ has a non-trivial type then the distance is
universally bounded. In addition, if $K \subset \RR^n$ is
quasi-convex then there exists a quasi-convex $T \subset \RR^n$
with a universally bounded isotropic constant and with a
universally bounded distance to $K$. }

\section{Introduction}

Let $K \subset \RR^n$ be a centrally symmetric (i.e. $K = -K$)
convex set with a non empty interior. Such sets are referred to
here as ``bodies''. We denote by $\langle \cdot, \cdot \rangle$
and $| \cdot |$ the standard scalar product and Euclidean norm in
$\RR^n$. We also write $D$ for the unit Euclidean ball and
$S^{n-1} =
\partial D$. The body $K$ has a unique linear image $\tilde{K}$
with $Vol(\tilde{K}) = 1$ such that
\begin{equation}
\int_{\tilde{K}} \langle x, \theta \rangle^2 dx \label{L_K_def}
\end{equation}
 does not depend on the choice of $\theta \in S^{n-1}$. We say
 that $\tilde{K}$ is an isotropic linear image of $K$ or that $\tilde{K}$
is in isotropic position. The quantity in (\ref{L_K_def}), for an
arbitrary $\theta \in S^{n-1}$, is usually referred to as $L_K^2$
or as the square of the isotropic constant of $K$. An equivalent
definition of $L_K$ which does not involve linear images is the
following:
\begin{equation}
n L_K^2 = \inf_T \int_K |Tx|^2 dx \label{L_K_def2}
\end{equation}
 where the infimum is over all
matrices $T$ such that $det(T) = 1$. For a comprehensive
discussion of the isotropic position and the isotropic constant we
refer the reader to \cite{MP}.

\medskip $L_K$ is an important linearly invariant
parameter associated with $K$. A major conjecture is whether there
exists a universal constant $c > 0$ such that $L_K < c$ for all
convex centrally symmetric bodies in all dimensions. A proof of
this conjecture will have various consequences. Among others (see
\cite{MP}) it will establish the fact that any body of volume one
has at least one $n-1$ dimensional section whose volume is greater
than some positive universal constant. This conjecture is known as
the slicing problem or the hyperplane conjecture. The best
estimate known to date is $L_K < c n^{1/4} \log n$ for $K \subset
\RR^n$ and is due to Bourgain \cite{Bou2} (see also the
presentation in \cite{D}). In addition, for large classes of
bodies the conjecture was positively verified (some examples of
references are \cite{Ba2}, \cite{Bou1}, \cite{KMP}, \cite{MP}).

\medskip In this note we deal with a known relaxation of this conjecture, which
we call the ``isomorphic slicing problem''. It was suggested to
the author by V. Milman. For two sets $K, T \subset \RR^n$, we
define their ``geometric distance'' as
$$ d_G(K,T) = \inf \left \{ ab ; \ \frac{1}{a} K \subset T \subset b K, \  a, b > 0
\right \}. $$ The Banach-Mazur distance of $K$ and $T$ is
$$ d_{BM}(K,T) = \inf \{ d_G(K, L(T)) \ ; \ L \ is \ a \ linear \ operator \}. $$
Let $K_n, T_n \subset \RR^n$ for $n=1,2,...$ be a sequence of
bodies such that $d_{BM}(K, T) < Const$ independent of the
dimension $n$. In this case we say that the families $\{K_n\}$ and
$\{T_n\}$ are uniformly isomorphic. Indeed, the norms that $K_n$
and $T_n$ are their unit balls are uniformly isomorphic. The
isomorphic slicing problem asks whether the slicing problem is
correct, at least up to a uniform isomorphism. Formally:

\begin{question} Does there exist constants $c_1, c_2 > 0$ such
that for any dimension $n$, for any body $K \subset \RR^n$, there
exists a body $T \subset \RR^n$ with $d_{BM}(K,T) < c_1$ and $L_T
< c_2$? \label{isomorphic_question}
\end{question}

In this note we answer this question affirmatively, up to a
logarithmic factor. The following is proved here:

\begin{theorem} For any centrally symmetric convex body $K \subset \RR^n$
there exists a centrally symmetric convex body $T \subset \RR^n$
with $d_{BM}(K, T) < c_1 \log n$ and
$$ L_T < c_2 $$
where $c_1, c_2 > 0$ are numerical constants. \label{main_result}
\end{theorem}

The $\log n$ factor in Theorem \ref{main_result} comes from the
use of the $l$-position and Pisier's estimate for the norm of the
Rademacher projection (see \cite{P}). Actually we prove, in the
notation of Theorem \ref{main_result}, that $d_{BM}(K, T) < c_1
M(K) M^*(K)$ (see definitions in Section \ref{function}).
Therefore we verify the validity of the isomorphic slicing
conjecture for bodies that have a linear image with bounded $M
M^*$. This is a large class of bodies, including all bodies with a
non trivial type. In addition, Proposition \ref{proj_of_bounded}
and Proposition \ref{exp_sqrt_n} provide other classes of bodies
for which Question \ref{isomorphic_question} has a positive
answer.

There exist some connections between the slicing problem and its
isomorphic versions. An example is provided in the following
lemma.

\begin{lemma}
Assume that there exist $c_1, c_2 > 0$ such that for any integer
$n$ and an isotropic body $K \subset \RR^n$ there exists an
isotropic body $T \subset \RR^n$ with $d_G(K, T) < c_1$ and $L_T <
c_2$. Then there exists $c_3 > 0$ such that for any integer $n$
and a body $K \subset \RR^n$, we have $L_K < c_3$.
\end{lemma}

\emph{Proof:} $L_T < c_2$, therefore $T$ is in $M$-position (as
observed by K. Ball, see definitions and proofs in \cite{MP}).
Since $d_G(K,T) < c_1$, also $K$ is in $M$-position. Using
Proposition 1.4 in \cite{BKM} we obtain a universal bound for the
isotropic constant. \hfill $\square$

\medskip A set $K \subset \RR^n$ is quasi-convex with constant $C
> 0$ if $conv(K) \subset C K$, where $conv$ denotes convex hull.
 For centrally symmetric quasi-convex sets, the isomorphic slicing problem
 has a positive answer. Formally, as is proved in Section
 \ref{quasi},

\begin{theorem}
For any $C > 1$, there exist $c_1,c_2 > 0$ with the following
property: If $K \subset \RR^n$ is centrally symmetric and
quasi-convex with constant $C$, then there exists a centrally
symmetric $T \subset \RR^n$ such that $d_{BM}(K, T) < c_1$ and
$L_T < c_2$. (Note that necesarily $T$ is $c_1 C$-quasi convex).
\label{quasi_theorem}
\end{theorem}

\smallskip Our proof has few consequences which are formulated
and proved in Section \ref{consequence}, among which is an
improvement of an estimate from \cite{BKM}. Throughout this paper
the letters $c, C, c^{\prime}, c_1,c_2, Const$ etc. denote some
positive numerical constants, whose value may differ in various
appearances. We ignore measurability issues as they are
non-essential to our discussion. All sets and functions used here
are assumed to be measurable.

\section{Log concave functions}
\label{concave}

In this section we collect a few facts regarding log-concave
functions, most of which are known and appear in \cite{Ba1} or
\cite{MP}, yet our versions are slightly different. $f: \RR^n
\rightarrow [0, \infty)$ is log-concave if $\log f$ is concave on
its support. $f$ is $s$-concave, for $s > 0$, if $f^{1/s}$ is
concave on its support. Any $s$-concave function is also
log-concave (see e.g. \cite{Bo}). Given a non-negative function
$f$ on $\RR^n$ we define for $x \in \RR^n$,
$$ \| x \|_f = \left( \int_0^{\infty} f \left(r x \right) r^{n+1} dr
\right)^{-1/{n+2}}. $$
We also define $K_f = \{ x \in \RR^n ; \| x \|_f \leq 1 \}$. The
following Busemann-type theorem appears in \cite{Ba1} (see also
\cite{MP}):

\begin{theorem} Let $f$ be an even log-concave function on
$\RR^n$. Then $K_f$ is convex and centrally symmetric and $\|
\cdot \|_f$ is a norm. \label{body_C}
\end{theorem}

In the sequel we make a repeated use of the following well known
facts. The first is that for any $1 \leq k \leq n$,
\begin{equation}
\left( \frac{n}{k} \right)^k \leq \left(\! \! \! \begin{array}{c}
n \\ k
\end{array} \! \! \! \right) < \left( e \frac{n}{k} \right)^k.
\label{gluskin}
\end{equation}
The second is that for any integers $a, b \geq 0$,
\begin{equation}
\label{beta}
 \int_0^1 s^a (1-s)^b ds = \frac{1}{(a+b+1) \left(\! \! \!
\begin{array}{c} a+b \\ a \end{array} \! \! \! \right) }.
\end{equation}

\begin{lemma} Let $f: \RR^n \rightarrow [0, \infty)$ be an even
function whose restriction to any straight line through the
origin is $s$-concave. Assume that $f(0) = 1$. If $s
> n$ then
$$ d_G(K_f, Supp(f)) < c \frac{s}{n} $$
where $c > 0$ is a numerical constant, and $Supp(f) = \{ x ; f(x)
> 0 \}$.
\label{distance_concave}
\end{lemma}

\emph{Proof:} Fix $\theta \in S^{n-1}$. Denote $M_\theta = \sup
 \{ r > 0; f(r \theta) > 0 \}$. Since $f|_{\theta \RR}$ is $s$-concave and $f(0)
 = 1$, for all $0 \leq r \leq M_{\theta}$,
 $$ f(r \theta) \geq \left( 1 - \frac{r}{M_{\theta}} \right)^s. $$
By the definition of $\| \theta \|_f$ and (\ref{beta}),
$$ \| \theta \|_f^{-(n+2)} \geq \int_0^{M_{\theta}} \left(1 -
\frac{r}{M_{\theta}} \right)^{s} r^{n+1} dr = \frac{M_{\theta}
^{n+2}}{(n+s+2) \left( \! \! \! \begin{array}{c} n+s+1 \\ n+1
\end{array} \! \! \! \right)}. $$
In addition, since $f|_{\theta \RR}$ is even, its maximum is $f(0)
= 1$ and
$$ \| \theta \|_f^{-(n+2)} \leq \int_0^{M_{\theta}} r^{n+1} dr =
\frac{1}{n+2} M_{\theta}^{n+2}. $$ Combining with the estimate
(\ref{gluskin}),
$$ \frac{(n+2)^{1/(n+2)}}{M_{\theta}} \leq \| \theta \|_f \leq
\frac{e (n+s+2)^{1/{n+2}} \left( \frac{n+s+1}{n+1}
\right)^{\frac{n+1}{n+2}}}{M_{\theta}}  $$ and because $s > n$,
$$ \forall \theta \in S^{n-1}, \ \frac{c_1}{M_{\theta}} < \| \theta \|_f
< \frac{c_2}{M_{\theta}} \frac{s}{n}
 \ \ \Rightarrow \ \ \frac{n}{c_2 s} Supp(f) \subset K_f \subset \frac{1}{c_1}
 Supp(f) $$
and the lemma is proved. \hfill $\square$

\medskip The isotropic constant and the isotropic position may be
defined for arbitrary measures or densities, rather than just for
convex bodies. Let $f:\RR^n \rightarrow [0, \infty)$ be an even
function with $\int_{\RR^n} f < \infty$. The entries of its
covariance matrix with respect to a fixed orthonormal basis
$\{e_1,..,e_n\}$ are defined as
$$ M_{i,j} = \frac{1}{\int_{\RR^n} f(x) dx} \int_{\RR^n} f(x) \langle x, e_i \rangle \langle x,
e_j \rangle dx. $$ We define $L_f = \left(
\frac{f(0)}{\int_{\RR^n} f} \right)^{\frac{1}{n}}
det(M)^{\frac{1}{2n}}$. One can verify that if $f$ is the
characteristic function of a body $K \subset \RR^n$, then $L_f =
L_K$. Our next lemma claims that if $f$ is log-concave, the body
$K_f$ shares the isotropic constant of the function $f$, up to a
universal constant. This fact appears in \cite{MP} and in
\cite{Ba1}, but not in a very explicit formulation. For
completeness we present a proof here.

\begin{lemma} Let $f$ be an even function on $\RR^n$ whose restriction to
any straight line through the origin is log-concave. Assume that
$\int_{\RR^n} f < \infty$. Then,
$$ c_1 L_f < L_{K_f} < c_2 L_f $$
where $c_1,c_2 > 0$ are universal constants.
 \label{L_equiv}
\end{lemma}

\emph{Proof:} Multiplying $f$ by a constant if necessary, we may
assume that $f(0) = 1$. Integrating in polar coordinates, for any
$y \in \RR^n$,
\begin{eqnarray*}
\lefteqn{\int_{K_f} \langle x, y \rangle^2 dx } \\
& = & \int_{S^{n-1}} \int_0^{1/\| \theta \|_f} \langle y, r \theta
\rangle^2 r^{n-1} dr d\theta = \frac{1}{n+2} \int_{S^{n-1}}
\langle y,
\theta \rangle^2 \frac{1}{\| \theta \|_f^{n+2}} d\theta \\
& = & \frac{1}{n+2} \int_0^{\infty} \int_{S^{n-1}} f(r \theta)
\langle y, \theta \rangle^2 r^{n+1} dr d\theta  = \frac{1}{n+2}
\int_{\RR^n} \langle x, y \rangle^2 f(x) dx
\end{eqnarray*}
where $d\theta$ is the induced surface area measure on
$S^{n-1}$. Denote by $M(f)$ and $M(K_f)$ the inertia matrices of
$f$ and the characteristic function of $K_f$, correspondingly.
Then $Vol(K_f) M(K_f) = \frac{1}{n+2} \left( \int_{\RR^n} f
\right) M(f)$. To compare the isotropic constants, we need to
estimate $\frac{\int f}{Vol(K_f)}$. Now,
\begin{equation}
Vol(K_f) = \frac{1}{n} \int_{S^{n-1}} \left( \int_0^{\infty} f
\left(r \theta \right) r^{n+1} dr \right)^{\frac{n}{n+2}} d\theta.
\label{Vol_polar}
\end{equation}
 We shall use the following one dimensional lemma, to be proved in the end of
this section (see also \cite{Ba1}, \cite{BKM} or \cite{MP}).

\begin{lemma}
Let $g:[0, \infty) \rightarrow [0, \infty)$ be a non-increasing
log-concave function with $g(0) = 1$ and $\int_0^{\infty} g(t)
t^{n-1} dt < \infty$. Then, for any integer $n \geq 1$,
$$ \frac{n^{\frac{n+2}{n}}}{n+2} \leq \frac{\int_0^{\infty}
g(t) t^{n+1} dt}{\left( \int_0^{\infty} g(t) t^{n-1} dt
\right)^{\frac{n+2}{n}}} \leq \frac{(n+1)!}{\left( (n-1)!
\right)^{\frac{n+2}{n}}}. $$ \label{one_dim}
\end{lemma}
(the left most inequality - which is more important for us -
is true without the log-concavity assumption).

Since $f$ is even and log-concave on any line through the origin,
it is non-increasing on any ray that starts in the origin. From
the left most inequality in Lemma \ref{one_dim}, for any $\theta
\in S^{n-1}$ (except for a set of measure zero where the integral
diverges),
$$ \int_0^{\infty} f \left(r \theta \right) r^{n+1}
dr \geq \frac{n^{\frac{n+2}{n}}}{n+2} \left( \int_0^{\infty} f
\left(r \theta \right) r^{n-1} dr \right)^{\frac{n+2}{n}} $$ and
according to (\ref{Vol_polar}),
$$ Vol(K_f) \geq \frac{1}{n} \frac{n^{\frac{n+2}{n}}}{n+2} \int_{S^{n-1}} \int_0^{\infty} f
\left(r \theta \right) r^{n-1} dr d\theta = \frac{n^{2/n}}{n+2}
\int_{\RR^n} f. $$ Since $M(K_f) = \frac{1}{n+2}
\frac{\int_{\RR^n} f}{Vol(K_f)} M(f)$,
$$ \frac{L_{K_f}^2}{L_f^2} = \frac{1}{n+2} \left(
\frac{\int_{\RR^n} f}{Vol(K_f)} \right)^{1 + \frac{2}{n}} \leq
\frac{1}{n+2} \left( \frac{n+2}{n^{2/n}} \right)^{\frac{n+2}{n}} <
c_2. $$ This finishes the proof of one part of the lemma. The
proof of the other inequality is similar. Using the right most
inequality in Lemma \ref{one_dim},
$$ \frac{L_{K_f}^2}{L_f^2} = \frac{1}{n+2} \left(
\frac{\int_{\RR^n} f}{Vol(K_f)} \right)^{1 + \frac{2}{n}} \geq
\frac{1}{n+2} \left( \frac{n \left( (n-1)!
\right)^{\frac{n+2}{n}}}{(n+1)!} \right)^{\frac{n+2}{n}} > c_1 $$
and the lemma is proved. \hfill $\square$

\medskip
\emph{Proof of Lemma \ref{one_dim}:} Start with the left-most
inequality. Define $A
> 0$ such that $ \int_0^{\infty} g(t) t^{n-1} dt = \int_0^A
t^{n-1} dt $. Then,
\begin{eqnarray*}
\lefteqn{\int_0^A (1 - g(t)) t^{n+1} dt - \int_A^{\infty} g(t) t^{n+1}dt } \\
& \leq & A^2 \left[ \int_0^A (1 - g(t)) t^{n-1} - \int_A^{\infty}
g(t) t^{n-1} dt \right] = 0.
\end{eqnarray*}
Since $\int_0^A t^{n+1} dt = \frac{n^{\frac{n+2}{n}}}{n+2} \left(
\int_0^A t^{n-1} dt \right)^{\frac{n+2}{n}}$, we get that
 $$ \int_0^{\infty} g(t) t^{n+1} dt \geq \int_0^A r^{n+1} dt = \frac{n^{\frac{n+2}{n}}}
 {n+2} \left( \int_0^{\infty} g(t) t^{n-1} dt \right)^{\frac{n+2}{n}}. $$
To obtain the other inequality we need to use the log-concavity of
the function. Define $B > 0$ such that $h(t) = e^{-B t}$ satisfies
$$ \int_0^{\infty} g(t) t^{n-1} dt = \int_0^{\infty} h(t) t^{n-1}
dt. $$ It is impossible that always $g < h$ or always $g > h$,
hence necessarily $t_0 = \inf \{ t > 0; h(t) \geq g(t) \}$ is
finite. $-\log g$ is convex and vanishes at zero, so $\tilde{g}(t)
= \frac{-\log g(t)}{t}$ is non-decreasing. Hence $(B -
\tilde{g}(t))(t - t_0) \geq 0$ or equivalently $(h(t)-g(t))(t-t_0)
\geq 0$ for all $t > 0$. Therefore,
\begin{eqnarray*}
\lefteqn{\int_0^{t_0} (g(t)-h(t)) t^{n+1} dt - \int_{t_0}^{\infty} (h(t) - g(t)) t^{n+1}dt } \\
& \leq & t_0^{2} \left[ \int_0^{t_0} (g(t) - h(t)) t^{n-1} -
\int_{t_0}^{\infty} (h(t) - g(t)) t^{n-1} dt \right] = 0.
\end{eqnarray*}
Since $\int_0^{\infty} e^{-t B} t^{n+1} dt = \frac{(n+1)!}{\left(
(n-1)! \right)^{\frac{n+2}{n}}} \left( \int_0^{\infty} e^{-t B}
t^{n-1} dt \right)^{\frac{n+2}{n}}$,
 $$ \int_0^{\infty} g(t) t^{n+1} dt \leq \int_0^{\infty} h(t) r^{n+1} dt =
 \frac{(n+1)!}{\left( (n-1)! \right)^{\frac{n+2}{n}}}   \left( \int_0^{\infty} g(t)
t^{n-1} dt \right)^{\frac{n+2}{n}}. $$ \hfill $\square$

\medskip Denote $L_n = \sup_{K \subset \RR^n} L_K$, where the
supremum is over all bodies in $\RR^n$. Define also $\tilde{L}_n =
\sup_{f:\RR^n \rightarrow [0, \infty)} L_f$ where the supremum is
over all log-concave even functions on $\RR^n$. Apriori, $L_n \leq
\tilde{L}_n$ since characteristic functions of bodies are
log-concave. Using the fact that convolution of log-concave
functions is again log-concave, and convolving the characteristic
function of $K$ with itself, we even obtain that $\sqrt{2} L_n
\leq \tilde{L}_n$. Less trivial is the fact that the opposite
inequality also holds, up to some constant (I was informed that
Corollary \ref{log_concave_as_bodies} also appears in K. Ball's
PhD thesis).

\begin{corollary} There exists $c > 0$ such that for any integer $n$,
$$ \tilde{L}_n < c L_n. $$
\label{log_concave_as_bodies}
\end{corollary}

\emph{Proof:} By lemma \ref{L_equiv}, for any log-concave even
$f:\RR^n \rightarrow [0, \infty)$,
$$ L_f < c L_{K_f} $$
and hence $\tilde{L}_n \leq c L_n$. \hfill $\square$

\section{Building a function on $K$}
\label{function}

Let $K \subset \RR^n$ be a body. In this section we find an
$\alpha n$-concave function $F$ supported on $K$ whose isotropic
constant is bounded. From Lemma \ref{L_equiv} it follows that
$L_{K_F} < Const$. According to Lemma \ref{body_C}, $K_F$ is a
convex body, and by Lemma \ref{distance_concave} we get that
$d_G(K, K_F) < c \alpha$. By obtaining good estimates on $\alpha$
Theorem \ref{main_result} would follow. Let $\| \cdot \|$ be the
norm that $K$ is its unit ball, and denote by $\sigma$ the unique
rotation invariant probability measure on $S^{n-1}$. The median of
$\| x \|$ on $S^{n-1}$ with respect to $\sigma$ is referred to as
$M^{\prime}(K)$. We abbreviate $M^{\prime} = M^{\prime}(K)$ and
define the following function on $K$:
 $$ f_K(x) = \inf \left \{ 0 \leq t \leq 1; x \in (1-t) \left
 [K \cap \frac{1}{M^{\prime}} D \right] + t K \right \}. $$
Then $f_K$ is a convex function which is zero on $K \cap
\frac{1}{M} D$. Define also
$$ M(K) = \int_{S^{n-1}} \| x \| d\sigma(x), \ \ M^*(K) = \int_{S^{n-1}} \| x \|_* d\sigma(x) $$
where $\| x \|_* = \sup_{y \in K} \langle x, y \rangle$ is the
dual norm. It is known (e.g. \cite{MS}) that $M(K)$ is comparable
to $M^{\prime}(K)$.

\begin{proposition}
Let $K \subset \RR^n$ be a body, and let $\alpha = c M(K) M^*(K)$. Then,
$$ \int_K \left(1 - f_K(x) \right)^{\alpha n} dx < 2 Vol
\left(K \cap \frac{1}{M^{\prime}} D \right) $$
where $c > 0$ is some numerical constant.
\label{building_f}
\end{proposition}

\emph{Proof:} We denote $F(x) = \left(1 - f(x)
\right)^{\alpha n}$. Then,
\begin{eqnarray*}
\lefteqn{ \int_K F(x) dx  = \int_0^1 Vol \{ x \in K ; F(x) \geq t
\} dt } \\ & = & \int_0^1 Vol \{ x \in K ; f(x) \leq 1 -
t^{\frac{1}{\alpha n}} \} dt
\end{eqnarray*}
and changing variables $s = 1 - t^{\frac{1}{\alpha n}}$ yields
$$ \int_K F(x) dx = \alpha n \int_0^1 (1-s)^{\alpha n - 1}
Vol \left( (1-s) \left[K \cap \frac{1}{M^{\prime}} D \right] + s K \right)
ds. $$ Expand the volume term into mixed volumes (see e.g.
\cite{Sch}):
$$ Vol \left( (1-s) \left[K \cap \frac{1}{M^{\prime}} D \right] + s K
\right) = \sum_{i=0}^n \left( \! \! \! \begin{array}{c} n \\ i
\end{array} \! \! \! \right) V_i s^i (1-s)^{n-i} $$
where $V_i = V(K, i; \left[K \cap \frac{1}{M^{\prime}} D \right], n-i)$. Therefore,
$$ \int_K F(x) dx = \alpha n \sum_{i=0}^n V_i \left( \! \! \! \begin{array}{c} n \\ i
\end{array} \! \! \! \right)  \int_0^1 s^i (1-s)^{ (\alpha + 1)
n - i - 1} ds $$ and by (\ref{beta}),
$$ \int_K F(x) dx = \frac{\alpha}{\alpha + 1} V_0 \sum_{i=0}^n
\frac{\left( \! \! \!
\begin{array}{ccc} n \\ i \end{array} \! \! \! \right)}
{\left( \! \! \! \begin{array}{ccc} (1+\alpha)n-1 \\ i \end{array}
\! \! \! \right)} \frac{V_i}{V_0}. $$ Using (\ref{gluskin}) we may write
\begin{equation}
\int_K F(x) dx = \frac{\alpha}{\alpha + 1} V_0 \left[1 +
\sum_{i=1}^n \left( c_{n,i} \frac{n}{(1 + \alpha)n - 1} \left
(\frac{V_i}{V_0}\right)^{1/i} \right)^i \right] \label{int_F}
\end{equation}
where $\frac{1}{e} \leq c_{n,i} \leq e$. By Alexandrov-Fenchel
inequalities $V_i^2 \geq V_{i-1} V_{i+1}$ for $i \geq 1$ (e.g.
\cite{Sch}). It follows that for $1 \leq i \leq j$,
\begin{equation}
 \left( \frac{V_i}{V_0} \right)^{1/i} \geq \left( \frac{V_j}{V_0}
\right)^{1/j}. \label{alexandrov}
\end{equation}
In particular, if $\alpha + 1 > 4 e \frac{V_1}{V_0}$, then by (\ref{alexandrov}),
$$ c_{n,i} \frac{n}{(1 + \alpha)n - 1} \left
(\frac{V_i}{V_0}\right)^{1/i} < \frac{2e}{1 + \alpha}
\frac{V_1}{V_0} \leq \frac{1}{2}. $$ Substituting in (\ref{int_F})
we obtain
$$ \int_K F(x) dx < V_0 \sum_{i=0}^n \frac{1}{2^i} < 2 V_0
= 2 Vol \left( K \cap \frac{1}{M^{\prime}} D \right). $$ It
remains only to show that our $\alpha = c M(K) M^*(K)$ is greater
than a constant times $\frac{V_1}{V_0}$. Since
$\frac{1}{M^{\prime}} D \cap K \subset \frac{1}{M^{\prime}} D$,
\begin{eqnarray*}
\lefteqn{ V_1 = V(K, 1 ; \left[K \cap \frac{1}{M^{\prime}} D \right], n-1) } \\
& \leq & V \left(K, 1; \frac{1}{M^{\prime}} D, n-1 \right) =
\frac{1}{(M^{\prime})^{n-1}} Vol(D) M^*(K)
\end{eqnarray*}
because $Vol(D) M^*(K) = V(K, 1; D, n-1)$ (see e.g. \cite{Sch}).
Regarding $V_0$, since $M^{\prime}$ is the median,
$$ \sigma \left( M^{\prime} K \cap S^{n-1} \right) \geq \frac{1}{2}
\ \ \ \Rightarrow \ \ \ Vol \left(K \cap \frac{1}{M^{\prime}} D
\right) \geq \frac{Vol \left(\frac{1}{M^{\prime}} D \right) }{2}.
$$ To conclude,
$$ \frac{V_1}{V_0} \leq \frac{1}{(M^{\prime})^{n-1}} Vol(D) M^*(K) \frac{2}{
\frac{1}{(M^{\prime})^n} Vol(D) } = 2 M^{\prime}(K) M^*(K). $$ The
median of a positive function is not larger than twice its
expectation, therefore $M^{\prime}(K) \leq 2 M(K)$, and we get
that for $\alpha = c M(K) M^*(K)$, it is true that $\alpha + 1 > 4
e \frac{V_1}{V_0}$ for a suitable numerical constant $c > 0$.
\hfill $\square$

\begin{corollary}
Let $K \subset \RR^n$ be a body, $\alpha = c M(K) M^*(K)$ and
denote $F(x) = \left( 1 - f_K(x) \right)^{\alpha n}$. Then,
$$ L_F < c^{\prime} $$
where $c$ is the constant from Proposition \ref{building_f} and
$c^{\prime} > 0$ is a numerical constant. \label{final_corr}
\end{corollary}

\emph{Proof:} Consider $F$ as a density on $K$, i.e. consider the
probability measure $\mu_F(A) = \frac{\int_A F(x) dx}{\int_K F(x)
dx}$. Since $F \equiv 1$ on $K \cap \frac{1}{M^{\prime}} D$, by
Proposition \ref{building_f},
$$ \mu \left( K \cap \frac{1}{M^{\prime}} D \right) > \frac{1}{2}. $$
In other words, the median of the Euclidean norm with respect to
$\mu$ is not larger than $\frac{1}{M^{\prime}}$. Since $F$ is
$\alpha n$-concave, by standard concentration inequalities for the
Euclidean norm with respect to log-concave measures (it follows,
e.g., from Theorem III.3 in \cite{MS}, due to Borell),
$$ \EE_{\mu} |x|^2 < \frac{c}{(M^{\prime})^2}. $$
Combining definition (\ref{L_K_def2}) and the fact that $L_F^2 =
\left(\frac{F(0)}{\int_K F} \right)^{\frac{2}{n}}
det(M_F)^{\frac{1}{n}}$ where $M_F$ is the covariance matrix, we
get that
$$ \left( \frac{\int_K F(x) dx }{F(0)} \right)^{\frac{2}{n}}
n L_F^2 \leq \EE_{\mu} |x|^2 < \frac{c}{(M^{\prime})^2}. $$ Since
$\int_K F(x) dx \geq Vol \left( \frac{1}{M^{\prime}} D \cap K
\right) \geq \frac{1}{2} Vol( \frac{1}{M^{\prime}} D )$ and $F(0)
= 1$, we obtain that $L_F^2 < \frac{c^{\prime}}{n Vol(D)^{2/n}} <
Const$. \hfill $\square$.

\medskip
\emph{Proof of Theorem \ref{main_result}:} We shall use the notion
of $l$-ellipsoid, and Pisier's estimate for $M(K) M^*(K)$. We
refer the reader to \cite{P} or \cite{MS} for definitions and
proofs. Let $K \subset \RR^n$ be a body. There exist a linear
image $\tilde{K}$ of $K$ such that its $l$-ellipsoid  is the
standard Euclidean ball. By Pisier's estimate,
$$ M^*(\tilde{K}) M(\tilde{K}) < c \log d_{BM}(K, D) < c^{\prime} \log n. $$
According to Corollary \ref{final_corr}, there exists an $\alpha
n$-concave function $F$ supported on $\tilde{K}$, with $\alpha = c
M(\tilde{K}) M^*(\tilde{K})$ and $L_F < c_1$. By Lemma
\ref{L_equiv} we get that $L_{K_F} < c_2$. From Lemma
\ref{distance_concave},
$$ d_{BM}(K, K_F) \leq d_G( \tilde{K}, K_F) < c \alpha < c^{\prime}
M(\tilde{K}) M^*(\tilde{K}) < C \log n. $$ This finishes the
proof. \hfill $\square$

\section{The quasi-convex case}
\label{quasi}

With an arbitrary body $K \subset \RR^n$ associated a special
ellipsoid, called an $M$-ellipsoid. An $M$-ellipsoid may be
defined by the following theorem (see \cite{M}, or chapter 7 in
the book \cite{P}):
\begin{theorem}
Let $K \subset \RR^n$ be a body. Then there exists an ellipsoid
$\E \subset \RR^n$ with $Vol( \E ) = Vol( K )$ such that
$$ Vol(K \cap \E)^{1/n} > c Vol(K)^{1/n} $$
where $c > 0$ is a numerical constant. We say that $\E$ is an
$M$-ellipsoid of $K$ (with constant $c$). \label{M_ellipsoid}
\end{theorem}

Let $K \subset \RR^n$ be a centrally symmetric quasi-convex body
with constant $C$ (in short ``a $C$-quasi-body''). Assume that
$Vol(K) = 1$ and that $conv(K)$ is in $M$-position. Moreover, we
may assume that $conv(K)$ is in a $1$-regular $M$-position, so
that $ conv(K) \subset Vol^{1/n}(conv(K)) n^2 D$ (e.g. \cite{P} or
\cite{GM}). Let us build the following function on $conv(K)$:
$$ F_K(x) = \left \{ \begin{array}{ccc} 1 & |x| \leq \sqrt{n} \\
\left(1 - \frac{|x| - \sqrt{n}}{M_x - \sqrt{n}} \right)^{\alpha n}
& |x| > \sqrt{n} \end{array}  \right. $$ for some $\alpha > 0$ to
be determined later, where
$$ M_x = \sup \left \{ r > 0; r \frac{x}{|x|} \in conv(K) \right \}. $$
$F_K$ is not log-concave, yet we may still consider the centrally
symmetric set $K_{F_K} \subset \RR^n$, defined in Section
\ref{concave}. Note that the restriction of $F_K$ to any straight
line through the origin is $\alpha n$-concave on its support,
hence it is possible to apply Lemma \ref{distance_concave} or
Lemma \ref{L_equiv}. We start with a one dimensional lemma.

\begin{lemma}
Let $0 < a < b$ and $\alpha > 1$ be such that $b > 2a \left( 1 +
\frac{\alpha}{e} \right)$. Let $n$ be a positive integer. Then,
$$ \int_a^b \left( 1 - \frac{t-a}{b-a} \right)^{\alpha n} t^n dt <
\left( \frac{c_1}{\alpha} \right)^n \int_a^b t^n dt $$ where $c_1
> 0$ is a numerical constant.
\label{one_dim2}
\end{lemma}

\emph{Proof:} Denote the integral on the left by $I$ and the
integral on the right by $J = \frac{1}{n+1} \left[b^{n+1} -
a^{n+1} \right]$. Changing variables $s = \frac{t-a}{b-a}$ we get
that
\begin{eqnarray*}
\lefteqn{I = (b-a) \int_0^1 (1 - s)^{\alpha n} \left( a + (b-a) s
\right)^n ds } \\ & = & (b-a) \sum_{i=0}^n \left( \! \! \!
 \begin{array}{c} n \\ i \end{array} \! \! \! \right) a^{n-i}
 (b-a)^i \int_0^1 (1-s)^{\alpha n}  s^i ds \\
\end{eqnarray*}
and using (\ref{beta}),
$$ I = (b-a) a^n \sum_{i=0}^n \frac{\left( \! \! \!
 \begin{array}{c} n \\ i \end{array} \! \! \! \right)}{(\alpha n +
 i + 1) \left( \! \! \!  \begin{array}{c} \alpha n + i \\ i \end{array} \!
\! \! \right)} \left(\frac{b-a}{a} \right)^i. $$ The estimate
(\ref{gluskin}) together with some trivial inequalities yields
that
$$ I \leq \frac{b-a}{\alpha n} a^n \sum_{i=0}^n \left(
\frac{e}{\alpha} \right)^i \left(\frac{b-a}{a} \right)^i =
\frac{b-a}{\alpha n} a^n \frac{q^{n+1} - 1}{q - 1} $$ where $q =
\frac{e(b-a)}{\alpha a}$. We assumed that $q \geq 2$, hence
$$ I \leq \frac{2}{en} (aq)^{n+1} = \frac{2}{en} \left( \frac{e}{\alpha} \right)^n (b-a)^{n+1}
< \left( \frac{c}{\alpha} \right)^n J. $$ \hfill $\square$

Next we show that for a suitable value of $\alpha$,
which is just a numerical constant, most of the mass of $F_K$ is not far from the origin.

\begin{lemma} For any $\alpha > 1$,
$$ \int_{\RR^n \setminus c_2 \sqrt{n} D} F_K(x) dx < \left( \frac{c_1}{\alpha}
\right)^{n-1} Vol(conv(K))
$$
where $c_1$ is the constant from Lemma \ref{one_dim2} and $c_2 = 2
\left( 1 + \frac{\alpha}{e} \right)
> 1$. \label{quasi_mass_close}
\end{lemma}

\emph{Proof:}  Note that,
$$ \int_{\RR^n \setminus \sqrt{n} D} F_K(x) dx =
\int_{S^{n-1}} \int_{\sqrt{n}}^{\max \{M_{\theta}, \sqrt{n} \}}
\left(1 - \frac{r - \sqrt{n}}{M_{\theta} - \sqrt{n}}
\right)^{\alpha n} r^{n-1} dr d\theta $$ where $d\theta$ is the
induced surface area measure over the sphere. Let $E = \{ \theta
\in S^{n-1}; M_{\theta} > c_2 \sqrt{n} \}$. By Lemma
\ref{one_dim2},
\begin{eqnarray*}
\lefteqn{\int_{\RR^n \setminus c_2 \sqrt{n} D} F_K(x) dx} \\
& < & \int_{E} \int_{\sqrt{n}}^{M_{\theta}} \left(1 - \frac{r -
\sqrt{n}}{M_{\theta} - \sqrt{n}} \right)^{\alpha n} r^{n-1} dr
d\theta \\
& < & \left( \frac{c_1}{\alpha} \right)^{n-1} \int_{E}
\int_{\sqrt{n}}^{M_{\theta}} r^{n-1} dr d\theta < \left(
\frac{c_1}{\alpha} \right)^{n-1} Vol(conv(K)).
\end{eqnarray*}
\hfill $\square$

\begin{lemma} Assume that $K \subset \RR^n$ is an $A$-quasi-body of volume one,
 and that $conv(K)$ is in a $1$-regular $M$-position
with constant $B$. Then for $\alpha = c_3(A,B)$,
$$ L_{F_K} < c_4(A,B) $$ where $c_3(A,B), c_4(A,B)$ depend solely on their parameters,
not on $K$ or $n$. \label{quasi_L_bounded}
\end{lemma}

\emph{Proof:} Note that $v^n := Vol(conv(K)) < A^n Vol(K) = A^n$.
By Brunn-Minkowski (e.g. \cite{Sch}), the function $f(t) =
Vol(conv(K) \cap t)$ is $n$-concave. Hence,
$$ Vol \left( conv(K) \cap \sqrt{n} D \right)^{1/n} \geq \frac{1}{c v}
Vol \left( conv(K) \cap c v \sqrt{n} D \right)^{1/n} > \frac{B
v}{c v} = \frac{B}{c} $$ where the constant $c > 1$ satisfies $c
\sqrt{n} Vol(D)^{1/n} > 1$, so that $Vol(c v \sqrt{n} D) >
Vol(conv(K))$. Let $\alpha > 2 c_1 c \frac{A}{B}$. By Lemma
\ref{quasi_mass_close},
$$ \int_{\RR^n \setminus c_2 \sqrt{n} D} F_K(x) dx < \left(
\frac{B}{2 c A} \right)^{n-1} Vol(conv(K)) $$ $$ < A \left(
\frac{B}{2 c} \right)^{n-1} Vol(K) < \frac{2 c A}{B}
\frac{1}{2^{n-1}} Vol \left( conv(K) \cap \sqrt{n} D \right) $$
for $c_2 = c^{\prime} \frac{A}{B}$. Define a measure by $\mu(E) =
\frac{\int_E F_K(x) dx}{\int_{\RR^n} F_K(x) dx}$. Since $F_K$
equals $1$ on $conv(K) \cap \sqrt{n} D$, we get that
$$ \mu(\RR^n \setminus c_2 \sqrt{n} D) < \frac{2 c A}{B}
\frac{1}{2^{n-1}}. $$ Since $conv(K) \subset Vol(conv(K))^{1/n}
n^2 D$,
$$ \EE_{\mu} |x|^2 < c_2^2 n + \frac{2 c A}{B}
\frac{1}{2^{n-1}} Vol(conv(K))^{1/n} n^2 < c^{\prime}
\frac{A^2}{B} n.
$$ Therefore $L_{F_K}^2 = L_{\mu}^2 < c^{\prime} \frac{A^2}{B} \left(
\frac{F_K(0)}{\int F_K} \right)^{\frac{2}{n}}$. Note that $F_K(0)
= 1$. Since $\int F_K \geq Vol(conv(K) \cap \sqrt{n} D)$, we
conclude that
$$ L_{F_K}^2 < c^{\prime} \frac{A^2}{B} \frac{c^2}{B^2} = c_4(A,B). $$
\hfill $\square$

\emph{Proof of Theorem \ref{quasi_theorem}}: Let $K \subset \RR^n$
be a $C$-quasi-body. Let $\tilde{K}$ be a linear image of $K$ such
that $Vol(\tilde{K}) = 1$ and $conv(\tilde{K})$ is in $1$-regular
$M$-position, with a universal constant $c > 0$. Consider the
function $F_{\tilde{K}}$ for $\alpha = c_3(C, c)$. By Lemma
\ref{distance_concave}, the body $T = \tilde{K}_{F_{\tilde{K}}}$
satisfies
$$ d_G(\tilde{K}, T) < c^{\prime} c_3 d_G(conv(\tilde{K}), T) < c^{\prime
\prime}(C)
$$
for some function $c^{\prime \prime}(C) > 0$. Also, by Lemma
\ref{L_equiv} and Lemma \ref{quasi_L_bounded},
$$ L_T < \tilde{c} L_{F_{\tilde{K}}} < \bar{c}(C) $$
for some $\bar{c}(C)$ a function of $C$. This completes the proof.
\hfill $\square$

\section{Consequences of the proof}
\label{consequence}

 Here we collect a few results which are byproducts of our
methods. Our first two propositions enrich the family of convex
bodies for which the ``isomorphic slicing problem'' has an
affirmative answer. In this section $Vol(T)$ denotes the volume of
a set $T \subset \RR^n$ relative to its affine hull.

\begin{lemma}
Let $K \subset \RR^n$ be an isotropic body of volume one, and let
$0 < \lambda < 1$ and $L_K < A$ for some $A > 1$. Then for any
subspace $E$ of dimension $\lambda n$,
$$ Vol(K \cap E)^{\frac{1}{n}} < c(A) $$
where $c(A)$ depends solely on $A$, and is independent of the body
$K$ and the dimension $n$. \label{L_section}
\end{lemma}

\emph{Proof:} Since $\EE_K |x|^2 < n A^2$, the median of $|x|$ on
$K$ is smaller than $2 \sqrt{n} A$. Denote $K^{\prime} = K \cap 2
\sqrt{n} A D$. Then $Vol(K^{\prime}) > \frac{1}{2}$. Also, given
any subspace $E \subset \RR^n$ of dimension $\lambda n$,
$$ Vol(K^{\prime} \cap E) \leq Vol(2 \sqrt{n} A D \cap E)
\leq \left(c \frac{A}{\sqrt{\lambda}} \right)^{\lambda n}. $$
Since $K^{\prime}$ is symmetric, $Vol(K^{\prime}) \leq
Vol(K^{\prime} \cap E) Vol(Proj_{E^{\perp}} K^{\prime})$, where
$E^{\perp}$ is the orthogonal complement of $E$ and
$Proj_{E^{\perp}}$ is the orthogonal projection onto $E^{\perp}$
in $\RR^n$. Therefore,
$$ Vol \left(Proj_{E^{\perp}} K \right) \geq Vol \left( Proj_{E^{\perp}}K^{\prime}\ \right)
\geq \frac{Vol(K^{\prime})}{Vol(K^{\prime} \cap E)} \geq \left(c
\frac{\sqrt{\lambda}}{A} \right)^{\lambda n}. $$ We denote the
polar body of $K$ by $K^{\circ} = \{ y \in \RR^n; \forall x \in K,
\langle x, y \rangle \leq 1 \}$. By Santal\'{o}'s inequality
\cite{Sa} and reverse Santal\'{o} \cite{BM} (recall that
projection and section are dual operations),
\begin{eqnarray}
\lefteqn{Vol(K \cap E) Vol \left(Proj_{E^{\perp}} K \right)} \\
& < & \left( \frac{c}{\lambda n} \right)^{\lambda n} \left(
\frac{c}{(1-\lambda) n} \right)^{(1-\lambda) n} \frac{1}{Vol
\left(Proj_E K^{\circ} \right) Vol(K^{\circ} \cap E^{\perp})} \nonumber \\
& < & \left( \frac{c^{\prime}}{n} \right)^n
\frac{1}{Vol(K^{\circ})} < \left( \frac{c^{\prime \prime}}{n}
\right)^n \frac{1}{Vol(D)^2} Vol(K) < \tilde{c}^n Vol(K)
\nonumber.
\end{eqnarray}
Hence,
$$ Vol(K \cap E)^{\frac{1}{n}} < \tilde{c} \frac{Vol(K)^{\frac{1}{n}}}{Vol \left(Proj_{E^{\perp}} K
\right)^{\frac{1}{n}}} < \tilde{c} \left( c
\frac{A}{\sqrt{\lambda}} \right)^{\lambda } < c^{\prime}
A^{\lambda}
$$ and the lemma is proven, with $c(A) = c A > c A^{\lambda}$.
\hfill $\square$

The following proposition states that the isomorphic slicing
conjecture holds for all projections to proportional dimension of
bodies with a bounded isotropic constant.

\begin{proposition} Let $K \subset \RR^n$ be a body with $L_K <
A$, and let $0 < \lambda < 1$. Then for any subspace $E$ of
dimension $\lambda n$, there exists a convex body $T \subset E$
such that
$$ d_{BM}(Proj_E(K), T) < c^{\prime}(\lambda), \ \ \ L_T < c(\lambda, A) $$
where $Proj_E$ is the orthogonal projection onto $E$ in $\RR^n$,
and $c^{\prime}(\lambda), c(\lambda, A)$ are functions independent
of $K$ or $n$. \label{proj_of_bounded}
\end{proposition}

\emph{Proof:} We may assume that $K$ is of volume one and in
isotropic position. For $x \in E$, define
$$ f(x) = Vol(K \cap [E^{\perp} + x]). $$
 For any $\theta_1, \theta_2 \in E$,
$$ \int_E \langle x, \theta_1 \rangle \langle x, \theta_2 \rangle
f(x) dx = \int_K \langle x, \theta_1 \rangle \langle x, \theta_2
\rangle dx. $$ Hence by Lemma \ref{L_section},
$$ L_f = (f(0))^{\frac{1}{\lambda n}} L_K <
 Vol(K \cap E^{\perp})^{\frac{1}{\lambda n}} A < c(A)^{\frac{1}{\lambda}} A =
c^{\prime}(\lambda, A) $$ Let $T = K_f$. By Lemma \ref{L_equiv} we
know that $L_T < \tilde{c} L_f < c^{\prime \prime}(\lambda, A)$.
Also, by Brunn-Minkowski (e.g. \cite{Sch}) $f$ is
$(1-\lambda)n$-concave. By Lemma \ref{distance_concave} $d_G(T,
Proj_E(K)) < c \frac{1 - \lambda}{\lambda}$. This finishes the
proof. \hfill $\square$.

\medskip
Our next proposition verifies the isomorphic slicing conjecture
under the condition that at least a small portion of $K$ (say, of
volume $e^{-\sqrt{n}}$) is located not too far from the origin.

\begin{proposition}
Let $K \subset \RR^n$ be a body of volume one, such that $K
\subset \beta n D$. Assume that $Vol(K \cap \gamma \sqrt{n} D) >
e^{-\delta \sqrt{n}}$. Then there exists a body $T \subset \RR^n$
such that
$$ d_{BM}(K, T) < c \left(1 + \frac{\beta \delta}{\gamma} \right), \ \ \ L_T < c^{\prime} \gamma $$
where $c, c^{\prime} > 0$ are numerical constants.
 \label{exp_sqrt_n}
\end{proposition}

\emph{Proof:} If $K \subset 2 \gamma \sqrt{n} D$ the proposition
is trivial since $L_K < c^{\prime} \gamma$. Assume the opposite,
and denote $C = K \cap 2 \gamma \sqrt{n} D$. Similarly to Section
\ref{function}, we define
$$ f(x) = \inf \{ 0 \leq t \leq 1 ; x \in (1-t) C + t K \} $$
and consider the density $F(x) = \left( 1 - f(x) \right)^{\alpha
n}$ on $K$. It is enough to show that $\frac{V(K, 1; C,
n-1))}{Vol(C)} < c \left(1 + \frac{\beta \delta}{\gamma} \right)$.
Indeed, in that case for $\alpha = c^{\prime} \left(1 +
\frac{\beta \delta}{\gamma} \right)$, as in Proposition
\ref{building_f}, we get that $\int_C F(x) dx> \frac{1}{2} \int_K
F(x) dx$ and the same argument as in Corollary \ref{final_corr}
shows that
$$ L_{K_F} < c^{\prime} \gamma, \ \ \
d_G(K_F, K) < c \left(1 + \frac{\beta \delta}{\gamma} \right). $$
Let us bound $\frac{V(K, 1; C, n-1)}{Vol(C)}$. Define $f(t) =
Vol(K \cap t D)$. According to our assumption, $\log f(\gamma
\sqrt{n})
> -\delta \sqrt{n}$ and $\log f(2\gamma \sqrt{n}) < 0$. We conclude that
there exists $ \gamma \sqrt{n} < t_0 < 2 \gamma \sqrt{n}$ with
$\left( \log f(t_0) \right)^{\prime} < \frac{\delta}{\gamma}$. By
Brunn-Minkowski inequality $\log f$ is concave, and $(\log
f)^{\prime}$ is decreasing. Therefore, for $t = 2 \gamma \sqrt{n}
\geq t_0$,
$$ (\log f(t))^{\prime} = \frac{Vol(K \cap t S^{n-1})}{Vol(K \cap t
D)} < \frac{\delta}{\gamma}. $$ For $x \in \partial C$, we denote
by $\nu_x$ the outer unit normal to $C$ at $x$, if it is unique.
Let $h_K(x) = \sup_{y \in K} \langle x, y \rangle$. Then (see
\cite{Sch}),
\begin{eqnarray*}
\lefteqn{V(K, 1; C, n-1) = \frac{1}{n} \int_{\partial C}
h_K(\nu_x) dx } \\ & = & \frac{1}{n} \int_{K \cap t S^{n-1}}
h_K(x) dx + \frac{1}{n} \int_{\partial C \setminus t S^{n-1}}
h_C(\nu_x) dx \\ & \leq & \frac{1}{n} \left( \frac{\delta}{\gamma}
Vol(C) \right) \beta n + Vol(C) = \left( 1 + \frac{\beta
\delta}{\gamma} \right) Vol(C)
\end{eqnarray*}
where we used the fact that $h_K \leq \beta n$ and that $Vol(C) =
\frac{1}{n} \int_{\partial C} h_C(\nu_x) dx$. This completes the
proof. \hfill $\square$

\medskip For $K \subset \RR^n$, the volume ratio of $K$ is defined
as,
$$ v.r.(K) = \sup_{\E \subset K} \left( \frac{Vol(K)}{Vol(\E)}
\right)^{\frac{1}{n}} $$ where the supremum is over all ellipsoids
contained in $K$. We denote
$$ L_n(a) = \sup \{ L_K \ ; \ K \subset \RR^n \ is \ a \ body, \ v.r.(K) \leq a \}. $$
In \cite{BKM} it is proved that for any $\delta > 0$,
\begin{equation}
L_n < c(\delta) \ L_n( v(\delta) )^{1 + \delta}.
\label{old_reduction}
\end{equation}
 where $c(\delta), v(\delta) \approx e^{\frac{c}{1-\delta}}$. Next, we
improve the dependence in (\ref{old_reduction}).

\begin{corollary} There exists $c_1, c_2 > 0$, such that for all
$n$,
$$ L_n < c_1 L_n(c_2). $$
\label{better_BKM}
\end{corollary}

\emph{Proof:} Our proof is a modification of the proof in
\cite{BKM}. As in that paper, let $K \subset \RR^n$ be a body such
that $L_K = L_n$ and $K$ is of volume one and in isotropic
position. As is proved in \cite{BKM}, there exists a subspace $F
\subset \RR^n$ of dimension $\lambda n = \left \lceil \frac{n}{4}
\right \rceil$ such that
$$ v.r.(Proj_E(K)) < c_1, \ \ \ \left(Vol(K \cap E^{\perp})\right)^{4/n} > c_2. $$
As in the proof of Proposition \ref{proj_of_bounded} here, for $x
\in E$, define
$$ f(x) = Vol(K \cap [E^{\perp} + x]). $$
Then $d_G(K_f, Proj_E(K)) < c \frac{1 - \lambda}{\lambda} < 3c$
and hence $v.r.(K_f) < 3c c_1$. Also,
$$ L_f = (f(0))^{\frac{1}{\lambda n}} L_K =
\left( Vol(K \cap E^{\perp}) \right)^{4/n} L_K > c^{\prime} L_K. $$
By Lemma \ref{L_equiv}, $L_{K_f} > \tilde{c} L_K = \tilde{c} L_n$.
To conclude,
$$ L_{\lambda n}(3 c c_1) \geq L_{K_f} > \tilde{c} L_n $$ and by
Proposition 1.3 in \cite{BKM} $L_n > c L_{\lambda n}$. Hence for
$m = \left \lceil \frac{n}{4} \right \rceil$ we have $L_m < c_1
L_m(c_2)$ and the corollary is proved. \hfill $\square$.

\medskip
\emph{Remark:} Currently, there is no good proven bound for $M(K)
M^*(K)$ in the non-symmetric case, and hence the central symmetry
assumption of the body is crucial to the proof of Theorem
\ref{main_result}. However, part of the statements in this paper
may be easily generalize to non-symmetric bodies. In particular,
Theorem \ref{quasi_theorem}, Corollary \ref{log_concave_as_bodies}, Proposition
\ref{proj_of_bounded}, Proposition \ref{exp_sqrt_n} and Corollary
\ref{better_BKM} also hold in the non-symmetric case.

\medskip
\emph{Acknowledgement:} I would like to thank Vitali Milman for
many excellent discussions regarding the slicing problems and
other problems in high dimensional geometry.


\begin{thebibliography}{99}
\label{sec:reference}

\bibitem[Ba1]{Ba1} K. M. Ball, Logarithmically concave functions and sections
of convex sets in $\RR^n$. Studia Math. {\bf 88} (1988) 69--84.

\bibitem[Ba2]{Ba2} K. M. Ball, Normed spaces with a
weak-Gordon-Lewis property, Proc. of Funct. Anal., University of
Texas and Austin (1987--1989), Lecture Notes in Math., vol. {\bf
1470}, Springer (1991) 36--47.

\bibitem[Bo]{Bo} C. Borell, Convex set functions in $d$-space.
Period. Math. Hungar. 6, no.2 (1975) 111--136.

\bibitem[Bou1]{Bou1} J. Bourgain, On high-dimensional maximal
functions associated to convex bodies.  Amer. J. Math.
108, no. 6  (1986),  1467--1476.

\bibitem[Bou2]{Bou2} J. Bourgain, On the distribution of polynomials on
high dimensional convex sets, Geometric aspects of functional
analysis (1989--90), Lecture Notes in Math., vol. {\bf 1469},
Springer Berlin (1991) 127--137.

\bibitem[BKM]{BKM} J. Bourgain, B. Klartag, V. Milman,
Symmetrization and isotropic constants of convex bodies, to appear
in Geometric aspects of functional analysis, Lecure Notes in Math.

\bibitem[BM]{BM} J. Bourgain, V. Milman, New volume ratio
properties for convex symmetric bodies in $\RR^n$, Invent. Math.
88 , no. 2 (1987) 319--340.

\bibitem[D]{D} S. Dar, Remarks on Bourgain's problem on slicing of
convex bodies, Geometric aspects of functional analysis, Operator
Theory: Advances and Applications, vol. {\bf 77} (1995) 61--66.

\bibitem[GM]{GM} A.A. Giannopoulos, V.D. Milman, Mean width and
diameter of proportional sections of a symmetric convex body, J.
Reine. angew. Math. {\bf 497} (1998) 113--139.

\bibitem[M]{M} V.D. Milman, In\'{e}galit\'{e} de Brunn-Minkowski inverse
et applications \`{a} le th\'{e}orie locale des espaces
norm\'{e}s. C.R. Acad. Sci. Paris, Ser. I {\bf 302} (1986) 25--28.

\bibitem[MP]{MP} V.D. Milman, A. Pajor, Isotropic position and
inertia ellipsoids and zonoids of the unit ball of a normed
$n$-dimensional space. Geometric aspects of functional analysis
(1987--88), Lecture Notes in Math., vol. {\bf 1376}, Springer
Berlin, (1989) 64--104.

\bibitem[MS]{MS} V.D. Milman, G. Schechtman, Asymptotic theory of
finite-dimensional normed spaces. Lecture Notes in Mathematics
{\bf. 1200}. Springer-Verlag, Berlin (1986).

\bibitem[KMP]{KMP} H. K\"{o}nig, M. Meyer, A. Pajor,
The isotropy constants of the Schatten classes are bounded.
Math. Ann. 312, no. 4 (1998) 773--783.

\bibitem[P]{P} G. Pisier, The volume of convex bodies and Banach
space geometry, Cambridge Tracts in Mathematics, Cambridge univ.
Press, vol. {\bf 94} (1997).

\bibitem[Sa]{Sa} L. A. Santal\'{o}, An affine invariant for
 convex bodies of $n$-dimensional space.
 (Spanish)  Portugaliae Math.  8,  (1949). 155--161.

\bibitem[Sch]{Sch} R. Schneider, Convex bodies: the Brunn-Minkowski theory.
Encyclopedia of Mathematics and its Applications {\bf 44},
Cambridge University Press, Cambridge (1993).

\end{thebibliography}
\end{document}